# HITTING TIMES FOR SPECIAL PATTERNS IN THE SYMMETRIC EXCLUSION PROCESS ON $Z^D$

By Amine Asselah and Paolo Dai Pra

*Université de Provence and Università di Padova*

We consider the symmetric exclusion process $\{\eta_t, t > 0\}$ on $\{0,1\}^{\mathbb{Z}^d}$. We fix a pattern $\mathcal{A} := \{\eta : \sum_\Lambda \eta(i) \geq k\}$, where $\Lambda$ is a finite subset of $\mathbb{Z}^d$ and $k$ is an integer, and we consider the problem of establishing sharp estimates for $\tau$, the hitting time of $\mathcal{A}$. We present a novel argument based on monotonicity which helps in some cases to obtain sharp tail asymptotics for $\tau$ in a simple way. Also, we characterize the trajectories $\{\eta_s, s \leq t\}$ conditioned on $\{\tau > t\}$.

**1. Introduction.** We consider the symmetric simple exclusion process (SSEP) on $\mathbb{Z}^d$, where particles are indistinguishable. The state space is $\Omega := \{\eta : \eta(i) \in \{0,1\} \text{ for } i \in \mathbb{Z}^d\}$ and a graphical construction of the process is as follows. To bonds of the cubic lattice $\mathbb{Z}^d$, we associate independent Poisson processes of intensity 1, at whose time realizations the contents of the corresponding adjacent sites are exchanged. We fix a local pattern $\mathcal{A} \subset \Omega$ that depends on $\{\eta(i) : i \in \Lambda\}$, where $\Lambda$ is a finite subset of $\mathbb{Z}^d$, and we consider the problem of establishing sharp estimates for the hitting time of $\mathcal{A}$, $\tau := \inf\{t : \eta_t \in \mathcal{A}\}$. For a physical motivation, see, for instance, [1]. The SSEP is a nonirreducible Markov process on an uncountable state space with the following special properties (enounced in greater generality than SSEP).

1. There is a partial order on the state space $\Omega$, say $\prec$.
2. The generator of the dynamics, $\mathcal{L}$, is monotone, that is, $e^{t\mathcal{L}}$ preserves increasing functions for any $t \geq 0$.
3. There is an invariant probability measure $\nu$ which satisfies the FKG inequality.
4. The pattern of interest, $\mathcal{A}$, is increasing, that is, $\xi \in \mathcal{A}$ and $\xi \prec \eta$ imply that $\eta \in \mathcal{A}$.

5. The dual $\mathcal{L}^*$ of $\mathcal{L}$ in $L^2(\nu)$ is monotone.

A simple consequence of properties 1–5 is the existence of a limit (see, e.g., [1], (2.7))

$$(1.1) \qquad \lambda = -\lim_{t\to\infty} \frac{1}{t}\log(P_\nu(\tau > t)).$$

However, to obtain estimates sharper than (1.1), in the context of particle systems satisfying 1–5 and in the case $\lambda$ is positive, is a more intricate matter. For this purpose, it is useful to study the regularity of generalized principal Dirichlet eigenfunctions, that is, probability measures $\mu$ with support in $\mathcal{A}^c$, satisfying, for every $\varphi$ in the domain of $\mathcal{L}$, denoted by $D(\mathcal{L})$, and $\varphi|_\mathcal{A} \equiv 0$,

$$(1.2) \qquad \int (\mathcal{L}\varphi + \lambda\varphi)\,d\mu = 0.$$

Measures satisfying (1.2) are also called quasi-stationary measures, since if we draw an initial configuration from any such measure, then, for any time $t > 0$, the law of $\eta_t$ conditioned on $\{\tau > t\}$ is time-ivariant. We denote by $T_t(\pi)$ the law of this conditioned process at time $t$ with initial probability measure $\pi$. We recall some works relevant to our context. First, some quasi-stationary measures are obtained as limits of linear combination of $\{T_t(\nu), t > 0\}$ (see Theorem 1 of [2] and Theorem 2.4 of [1]). Assume that such a limit $\mu$ is absolutely continuous with respect to $\nu$, and call its density $u := d\mu/d\nu$. When $\mathcal{L}^*$ generates a Markov process, let $\mu^*$ be its corresponding quasi-stationary measure and assume it has a density $u^* := d\mu^*/d\nu$. In [1], Corollary 2.8 and its proof, we have the following general fact.

FACT 1.1. Assume that $\lambda$ given in (1.1) is positive and $u, u^* \in L^p(\nu)$ for $p > 2$. Then, for any $t \geq 0$,

$$(1.3) \qquad \exp(-H(\tilde{\nu}, \nu)) \leq \frac{P_\nu(\tau > t)}{\exp(-\lambda t)} \leq 1$$

with

$$d\tilde{\nu} = \frac{uu^*\,d\nu}{\int uu^*\,d\nu} \quad \text{and} \quad H(\tilde{\nu}, \nu) = \int \log\left(\frac{d\tilde{\nu}}{d\nu}\right) d\tilde{\nu} < \infty.$$

In the symmetric case, the results are stronger (see [2] or [3], Corollary 2.5).

FACT 1.2. If $\mathcal{L}$ is a self-adjoint Markov generator on $L^2(\nu)$, $\lambda > 0$ and $u \in L^2(\nu)$, then

$$(1.4)\ \lim_{t\to\infty} \frac{P_\nu(\tau > t)}{\exp(-\lambda t)} = \frac{(\int u\,d\nu)^2}{\int u^2\,d\nu} \qquad \text{with } \lambda = \inf_{f \in D(\mathcal{L})}\left\{\frac{-\int f\mathcal{L}f\,d\nu}{\int f^2\,d\nu} : f|_\mathcal{A} = 0\right\}.$$

Now, a key step in the proof of the regularity of quasi-stationary measures is to obtain uniform estimates for $\{T_t(\nu), t > 0\}$. In other words, we look for measures $\underline{\nu}$ and $\overline{\nu}$ such that, for any $t > 0$,

$$\underline{\nu} \prec T_t(\nu) \prec \overline{\nu} \qquad \left(\mu \prec \nu \text{ means that } \int f\, d\mu \leq \int f\, d\nu \text{ for all increasing } f\right)$$

and with $d\underline{\nu}/d\nu$ and $d\overline{\nu}/d\nu$ regular enough (see, e.g., [2] and [3]). In Section 2, we present a simple method to obtain such uniform stochastic bounds. Roughly, the main idea is to bound the principal eigenfunction $u$—which satisfies on $\mathcal{A}^c$ that $\mathcal{L}(u)/u$ is constant—by a *simple* function $\psi$ on which we impose a weaker assumption, namely that $\mathcal{L}(\psi)/\psi$ is increasing on $\mathcal{A}^c$. We first apply this method, in Section 3.3, to the SSEP on $\mathbb{Z}^d$ and the pattern $\mathcal{A}_1 := \{\eta : \eta(0) = 1\}$. In this context, $\xi \prec \eta$ when $\xi(i) \leq \eta(i) \ \forall i \in \mathbb{Z}^d$. Also, we recall that the SSEP has a one-parameter family of ergodic invariant measures $\{\nu_\rho : \rho \in [0,1]\}$, where $\nu_\rho$ is a product of Bernoulli measures of density $\rho$.

Thus, our first application is a key result of [3].

PROPOSITION 1.3. *Consider the SSEP in dimension $d \geq 5$, with pattern $\mathcal{A}_1$. For any density $\rho \in\, ]0,1[$, there is a sequence $\{\alpha_i, i \in \mathbb{Z}^d\}$ and a probability density $\psi$ with $\alpha_i \leq \rho$ for all $i \in \mathbb{Z}^d$,*

(1.5)
$$\sum_{i \in \mathbb{Z}^d \setminus \{0\}} \left(1 - \frac{\alpha_i}{\rho}\right)^2 < \infty,$$

$$\psi(\eta) := \frac{1}{Z}(1 - \eta(0)) \prod_{i \in \mathbb{Z}^d \setminus \{0\}} \left(\frac{\alpha_i}{1 - \alpha_i} \frac{1 - \rho}{\rho}\right)^{\eta(i)}$$

*($Z$ is a normalizing constant) such that if $d\nu_\alpha := \psi\, d\nu_\rho$, then for any $t > 0$,*

(1.6) $$\nu_\alpha \prec T_t(\nu_\rho) \prec \nu_\rho.$$

A corollary of Proposition 1.3 (see [3], Lemma 2.3) is the existence of $\mu_\rho := \lim_{t \to +\infty} T_t(\nu_\rho)$ as a strong limit in $L^2(\nu_\rho)$, that is, $dT_t(\nu_\rho)/d\nu_\rho$ converges in $L^2(\nu_\rho)$ to $d\mu_\rho/d\nu_\rho$. This $\mu_\rho$ is a quasi-stationary measure and is referred to as a *Yaglom limit*.

As a second illustration, we treat, in Section 3.4, the pattern $\mathcal{A}_2 := \{\eta : \eta(0) = \eta(0') = 1\}$, where $0'$ is a neighbor of the origin 0. However, for technical reasons, we need to have an intensity rate between 0 and $0'$ larger than $2d - 1$.

PROPOSITION 1.4. *Let $T_t^\beta(\nu_\rho)$ be the law at time $t$ of the SSEP modified by letting $\beta$ be the intensity rate between $(0, 0')$ and conditioned on $\{\tau > t\}$ with initial measure $\nu_\rho$. If the dimension $d \geq 5$ and $\beta \geq 2d - 1$, then stochastic estimates of type (1.6) hold.*

REMARK 1.5. To explain the reason for speeding up the intensity of bond $(0,0')$, we need to unravel a key technical assumption. The above mentioned function $\psi$, which mimics the Dirichlet eigenfunction, is associated with a Markov process that never enters $\mathcal{A}$ and has a formal generator

$$\mathcal{L}_\psi(\varphi) = \frac{\mathcal{L}(\psi\varphi) - \varphi\mathcal{L}(\psi)}{\psi}.$$

A handy assumption on $\mathcal{L}_\psi$ is that it is monotone. This fails to be the case for SSEP with $\mathcal{A}_2 = \{\eta : \eta(0) = \eta(0') = 1\}$. In other words, there is no coupling of two trajectories $(\eta_\cdot, \zeta_\cdot)$ governed by $\mathcal{L}_\psi$, with $\zeta_0 \prec \eta_0$, where the order is preserved in time. Indeed, consider $\zeta \prec \eta$ with $\eta(0'') = \eta(0') = 1$, where $0''$ is a neighbor of 0 different from $0'$, and $\zeta(0'') = 1 = 1 - \zeta(0')$. For the configuration $\eta$, the rate intensity associated with $(0,0'')$ is null, whereas it is positive in the configuration $\zeta$. Thus, if the first time realization of the Poisson process associated with $(0,0'')$ in $\zeta$ occurs before realizations of the processes associated with the other bonds adjacent to $0''$, then the order is destroyed. We show that speeding up the intensity of the process associated with $(0,0')$ enables us to build a monotone coupling.

Our method can also be used to prove regularity of invariant measures. Thus, our final application, in Section 3.5, is to study the regularity of invariant measures for the symmetric exclusion dynamics with birth and death of particles at the origin. For simplicity, we consider the process where the neighbors of the origin can die with positive rate $a$ and be born with positive rate $b$. The invariant measures were studied in [6]. We obtain here a new characterization.

PROPOSITION 1.6. *When $d \geq 5$, there is a stationary measure $\mu_\rho^{ab}$, for any $\rho \in\, ]0,1[$ such that*

$$\bigotimes_{i \neq 0} \nu_{\alpha_i} \prec \mu_\rho^{ab} \prec \bigotimes_{i \neq 0} \nu_{\tilde{\alpha}_i} \quad \text{and}$$

(1.7)
$$1 + C_{ab}\mathbb{P}_i(H_0 < \infty) = \frac{\tilde{\alpha}_i}{\rho}\frac{1-\rho}{1-\tilde{\alpha}_i} = \frac{\rho}{\alpha_i}\frac{1-\alpha_i}{1-\rho},$$

*where $\mathbb{P}_i(H_0 < \infty)$ is the probability that a symmetric random walk starting at site $i$ hits the origin, $C_{ab}$ is a positive constant depending on $a$ and $b$, and $\bigotimes_{i \neq 0} \nu_{\alpha_i}$ denotes a product Bernoulli measure of density $\alpha_i$ at site $i$ of $\mathbb{Z}^d \setminus \{0\}$.*

REMARK 1.7. This implies by the arguments of [2] that $\mu_\rho^{ab}$ is equivalent to $\nu_\rho$ and that $d\mu_\rho^{ab}/d\nu_\rho$ is in $L^p(\nu_\rho)$ for any integer $p$ when $d \geq 5$.

The problems we consider in Sections 4 and 5 are inspired by works on conditional Brownian motion (see, e.g., [4], Theorems 1 and 2, [12], Theorem 3 and [10]). We assume the following hypotheses.

HYPOTHESES ($\mathcal{H}$). *The generator $\mathcal{L}$ is self-adjoint in $L^2(\nu)$. The Yaglom limit $\mu := \lim_{t \to +\infty} T_t(\nu)$ exists with a corresponding $\lambda > 0$ for which* (1.1) *holds. Moreover, $u := d\mu/d\nu \in L^2(\nu)$, $u$ is a simple eigenfunction for $\lambda$, $u$ is positive $\nu$-a.s. and*

$$\text{(1.8)} \qquad \frac{dT_t(\nu)}{d\nu} \xrightarrow{L^2(\nu)} u.$$

Hypotheses ($\mathcal{H}$) were proved in [3] for SSEP in dimension $d \geq 5$ with $\mathcal{A}_1$. For the pattern $\mathcal{A}_2$, although the convergence in (1.8) is a corollary of Proposition 1.4, the uniqueness of $u$ in $L^2(\nu_\rho)$ is open.

PROPOSITION 1.8. *Assume ($\mathcal{H}$).*

(i) *For every $f, g \in L^2(\nu)$, we have*

$$\text{(1.9)} \qquad \lim_{t \to \infty} \frac{E_\nu[f(\eta_0) g(\eta_t) \mathbb{1}_{\{\tau > t\}}]}{P_\nu(\tau > t)} = \int f \, d\mu \int g \, d\mu.$$

(ii) *For any measure $\pi$ with $d\pi/d\nu \in L^2(\nu)$, we have the weak-$L^2(\nu)$ convergence*

$$T_t(\pi) \xrightarrow{t \to \infty} \mu.$$

Finally, let $d\hat{\mu} = u^2 \, d\nu / \int u^2 \, d\nu$ and let $\{P^u_\eta, \eta \in \Omega\}$ be the law of the Markov process, reversible in $L^2(\hat{\mu})$, formally generated on $\mathcal{A}^c$ by

$$\mathcal{L}_u \varphi = \frac{\mathcal{L}(u\varphi) - \varphi \mathcal{L}(u)}{u} \qquad \text{(see definition in Section 5).}$$

We have the following characterization of trajectories in $\{\tau > t\}$.

PROPOSITION 1.9. *Assume ($\mathcal{H}$). Let $t \mapsto a_t$ be an increasing positive function such that $\lim_{t \to \infty} a_t = \lim_{t \to \infty}(t - a_t) = \infty$. For any $r > 0$, the law of $\{\eta_{a_t+s}, s \in [0, r]\}$, conditioned on $\{\tau > t\}$ with initial measure $\nu$, converges to the restriction to the time interval $[0, r]$ of $\int P^u_\eta \, d\hat{\mu}(\eta)$ (convergence in the topology induced by duality against bounded measurable functions).*

**2. The monotone method.** We consider a finite state space $\mathcal{X}$ with partial order $\prec$. We recall that a dynamics is monotone when its evolution semigroup preserves increasing functions or, equivalently, when there is a coupling of two paths $(\eta_t, \zeta_t)$ such that if $\eta_0 \prec \zeta_0$, then $P(\eta_t \prec \zeta_t \, \forall t \geq 0) = 1$.

Let $\{P_\eta(\cdot), \eta \in \mathcal{X}\}$ be a Markov process on $\mathcal{X}$ and let $\mathcal{L}$ be the corresponding infinitesimal generator.

LEMMA 2.1. *Let $\mathcal{A} \subset \mathcal{X}$ and $\tau = \inf\{t : \eta_t \in \mathcal{A}\}$. Assume that there is a function $\psi$ satisfying* (i) *$\psi$ is positive on $\mathcal{A}^c$ and $\psi|_\mathcal{A} = 0$,* (ii) *$\psi$ is decreasing on $\mathcal{A}^c$,* (iii) *$\mathcal{L}(\psi)/\psi$ is increasing on $\mathcal{A}^c$,* (iv) *the following Markov generator on $\mathcal{A}^c$ is monotone*:

$$\mathcal{L}_\psi(\varphi) := \frac{\mathcal{L}(\psi\varphi) - \varphi\mathcal{L}(\psi)}{\psi}. \tag{2.1}$$

*Then $\eta \mapsto P_\eta(\tau > t)/\psi(\eta)$ is increasing.*

PROOF. If $\{c(a,b), a,b \in \mathcal{X}\}$ are the rates associated with $\mathcal{L}$, then after a simple computation,

$$\forall a \in \mathcal{A}^c \quad \mathcal{L}_\psi f(a) = \sum_{b \in \mathcal{X}} c(a,b) \frac{\psi(b)}{\psi(a)}(f(b) - f(a))$$

$$= \sum_{b \notin \mathcal{A}} c(a,b) \frac{\psi(b)}{\psi(a)}(f(b) - f(a)). \tag{2.2}$$

Thus, $\mathcal{L}_\psi$ generates a Markov process on $\mathcal{A}^c$. By definition, for any $\varphi|_\mathcal{A} \equiv 0$,

$$\frac{\mathbb{1}_{\mathcal{A}^c}\mathcal{L}(\psi\varphi)}{\psi} = \mathbb{1}_{\mathcal{A}^c}\left(\mathcal{L}_\psi(\varphi) + \frac{\mathcal{L}(\psi)}{\psi}\varphi\right)$$

$$\implies \frac{\exp(t\mathbb{1}_{\mathcal{A}^c}\mathcal{L})(\psi\varphi)}{\psi} = \exp\left(t\mathbb{1}_{\mathcal{A}^c}\left(\mathcal{L}_\psi + \frac{\mathcal{L}\psi}{\psi}\right)\right)\varphi. \tag{2.3}$$

If $\{P_\eta^\psi(\cdot), \eta \in \mathcal{A}^c\}$ corresponds to $\mathcal{L}_\psi$, then (2.3) and the Feynmann–Kac formula give, for $\eta \notin \mathcal{A}$,

$$\frac{\int \varphi(\eta_t)\psi(\eta_t)\mathbb{1}_{\tau>t}\,dP_\eta}{\psi(\eta)} = \int \varphi(\eta_t)\mathbb{1}_{\tau>t}\exp\left(\int_0^t \frac{\mathcal{L}\psi}{\psi}(\eta_s)\,ds\right)dP_\eta^\psi$$

$$= \int \varphi(\eta_t)\exp\left(\int_0^t \frac{\mathcal{L}\psi}{\psi}(\eta_s)\,ds\right)dP_\eta^\psi. \tag{2.4}$$

Thus, for $\varphi = 1/\psi$,

$$\frac{P_\eta(\tau > t)}{\psi(\eta)} = \int \frac{1}{\psi(\eta_t)}\exp\left(\int_0^t \frac{\mathcal{L}\psi}{\psi}(\eta_s)\,ds\right)dP_\eta^\psi. \tag{2.5}$$

From (2.5), the lemma is proved using (ii)–(iv). $\square$

We state a related result. Assume that $\mathcal{L}$ generates an irreducible Markov process on $\mathcal{X}$ and let $\nu$ be a positive probability on $\mathcal{X}$. Denote by $\mathcal{L}^*$ the dual of $\mathcal{L}$ in $L^2(\nu)$. Note that $\mathcal{L}^*$ is not necessarily a Markov generator [since $\mathcal{L}^*(1) \neq 0$] and that by the Perron–Frobenius theorem (see, e.g., [11], Theorem 9.34), there is $u > 0$ with $\mathcal{L}^*(u) = 0$.

LEMMA 2.2. *Assume there is a function $\psi$ satisfying* (i) *$\psi$ is positive*, (ii) *$\mathcal{L}^*(\psi)/\psi$ is increasing and* (iii) *the following Markov generator is monotone*:

$$\mathcal{L}_\psi(\varphi) := \frac{\mathcal{L}^*(\psi\varphi) - \varphi\mathcal{L}^*(\psi)}{\psi}. \tag{2.6}$$

*Then $u/\psi$ is increasing. Similarly, if we assume $\psi'$ positive, $\mathcal{L}^*(\psi')/\psi'$ decreasing and $\mathcal{L}_{\psi'}$ monotone, then we obtain that $u/\psi'$ is decreasing.*

PROOF. We call $\varphi = u/\psi$ and look for the equation solved by $\varphi$:

$$\frac{\mathcal{L}^*(\varphi\psi)}{\psi} = 0 \implies \mathcal{L}_\psi(\varphi) + \frac{\mathcal{L}^*\psi}{\psi}\varphi = 0. \tag{2.7}$$

Note also that $\varphi$ is the principal eigenfunction of $\mathcal{L}_\psi + \mathcal{L}^*(\psi)/\psi$. By the Perron–Frobenius theorem and the Feynmann–Kac formula,

$$\begin{aligned}\varphi(\eta) &= \lim_{t\to\infty} \exp\left(t\left(\mathcal{L}_\psi + \frac{\mathcal{L}^*\psi}{\psi}\right)\right)1(\eta) \\ &= \lim_{t\to\infty} \int \exp\left(\int_0^t \frac{\mathcal{L}^*\psi}{\psi}(\eta_s)\,ds\right) dP_\eta^\psi.\end{aligned} \tag{2.8}$$

By hypotheses (ii) and (iii), we obtain that $\varphi$ is increasing.

With the same reasoning,

$$\frac{u}{\psi'} = \lim_{t\to\infty} \int \exp\left(\int_0^t \frac{\mathcal{L}^*\psi'}{\psi'}(\eta_s)\,ds\right) dP_\eta^{\psi'}$$

is decreasing since $\mathcal{L}^*(\psi')/\psi'$ is decreasing. □

**3. Three applications.** We consider three applications of the lemmas of Section 2. In Section 3.1, we introduce three particle systems. In proving Propositions 1.3, 1.4 and 1.6, the first step, carried out in Section 3.2, is to approximate these particle systems by finite-dimensional irreducible dynamics. The second step is to verify the hypotheses of Lemma 2.1 or 2.2 in each of our three cases. This is carried out, respectively, in Sections 3.3–3.5.

3.1. *Models.* First, we consider SSEP on $\Omega$ with the generator acting on local functions as

$$\mathcal{L}_{se}\varphi(\eta) = \sum_{i\in\mathbb{Z}^d}\sum_{j\sim i}(\varphi(T^{i,j}\eta) - \varphi(\eta)),$$

where $i \sim j$ means that $|i_1 - j_1| + \cdots + |i_d - j_d| = 1$, and

$$T^{i,j}\eta(j) = \eta(i), \qquad T^{i,j}\eta(i) = \eta(j) \quad \text{and for } k \neq i,j,\ T^{i,j}\eta(k) = \eta(k).$$

It is well known ([9], Theorem 3.9 and Example 3.1(d)) that $\mathcal{L}_{se}$ generates a Feller process and that the following set is a core of continuous functions:

$$\mathcal{D} := \left\{ \varphi : \sum_{i \in \mathbb{Z}^d} \nabla_i(\varphi) < \infty \right\}$$

where $\nabla_i(\varphi) = \sup\{|\varphi(\eta) - \varphi(\xi)| : \eta(j) = \xi(j) \ \forall j \neq i\}$.

It is also well known that for any $\rho \in [0,1]$, $\mathcal{L}_{se}$ extends to a self-adjoint operator on $L^2(\nu_\rho)$ (see, e.g., Section 2 of [13]).

Second, to treat $\mathcal{A}_2 := \{\eta : \eta(0) = \eta(0') = 1\}$, where $0'$ is a given neighbor of 0, we need to modify the intensity between the bond $b := (0, 0')$. Thus, we consider the generator

$$(3.1) \qquad \mathcal{L}_\beta \varphi = \mathcal{L}_{se} \varphi + (\beta - 1)(\varphi \circ T^b - \varphi) \qquad \text{with } \beta > 2d - 1.$$

Note that $\mathcal{L}_\beta$ is still self-adjoint in $L^2(\nu_\rho)$, for any $\rho \in \,]0,1[$.

Finally, we consider SSEP with birth and death of particles at neighbors of the origin. Thus, the state space is $\Omega^* := \{\eta(i) \in \{0,1\}, i \in \mathbb{Z}^d \setminus \{0\}\}$ and if $\mathcal{N}_0 := \{i \in \mathbb{Z}^d : i \sim 0\}$, then the generator $\mathcal{L}_{ab}$ reads as

$$(3.2) \quad \begin{aligned} \mathcal{L}_{ab}\varphi(\eta) &= \sum_{e \notin \mathcal{N}_0 \times \{0\}} (\varphi(T^e \eta) - \varphi(\eta)) \\ &\quad + \sum_{k \sim 0} (a\eta(k) + b(1 - \eta(k)))(\varphi(\sigma_k \eta) - \varphi(\eta)), \end{aligned}$$

where $\sigma_k$ is the spin flip at site $k$, $\sigma_k \eta(k) = 1 - \eta(k)$ and $\sigma_k \eta(j) = \eta(j)$ for $j \neq k$.

3.2. *Approximation by irreducible dynamics.* Let $\Lambda_n := [-n, n]^d$ and $\mathcal{A} \subset \Omega_n := \{0,1\}^{\Lambda_n}$. For a subset $U \subset \mathbb{Z}^d$, we denote by $\mathcal{F}_U$ the $\sigma$-field generated by $\{\eta(i), i \in U\}$. We set, for $\varphi$ on $\Omega_n$,

$$(3.3) \qquad \mathcal{L}_{se}^{n,\rho}\varphi := E_{\nu_\rho}[\mathcal{L}_{se}\varphi | \mathcal{F}_{\Lambda_n}] \quad \text{and} \quad \mathcal{L}_\beta^{n,\rho}\varphi := E_{\nu_\rho}[\mathcal{L}_\beta \varphi | \mathcal{F}_{\Lambda_n}].$$

For $\varphi$ on $\Omega_n^* := \{0,1\}^{\Lambda_n \setminus \{0\}}$, we set

$$(3.4) \qquad \mathcal{L}_{ab}^{n,\rho}\varphi := E_{\nu_\rho}[\mathcal{L}_{ab}\varphi | \mathcal{F}_{\Lambda_n \setminus \{0\}}].$$

An easy computation gives

$$(3.5) \quad \begin{aligned} \mathcal{L}_{se}^{n,\rho}\varphi(\eta) &= \sum_{\substack{i \sim j \\ i,j \in \Lambda_n}} (\varphi(T^{i,j}\eta) - \varphi(\eta)) \\ &\quad + \sum_{i \in \partial \Lambda_n} n(i) \sqrt{\frac{d\sigma_i \nu_\rho}{d\nu_\rho}(\eta)} (\varphi(\sigma_i \eta) - \varphi(\eta)), \end{aligned}$$

where $\partial \Lambda_n := \{i \in \Lambda_n : \exists j \notin \Lambda_n \text{ with } j \sim i\}$ and $n(i) = |\{j \notin \Lambda_n : j \sim i\}|$. A similar formula holds for $\mathcal{L}^{n,\rho}_\beta$. It follows easily from their definition that $\mathcal{L}^{n,\rho}_{se}$ and $\mathcal{L}^{n,\rho}_\beta$ are $(\nu_\rho|_{\Lambda_n})$-reversible on $\Omega_n$. We state next the irreducibility property, although the immediate proof is omitted.

LEMMA 3.1. *The generator $\mathcal{L}^{n,\rho}_{ab}$ is irreducible on $\Omega_n^*$.*

The dual of $\mathcal{L}^{n,\rho}_{ab}$ in $L^2(\Omega_n^*, \nu_\rho)$ is obtained after the simple computation

$$(3.6) \quad (\mathcal{L}^{n,\rho}_{ab})^* f(\eta) = \mathcal{L}_0 f(\eta) + \sum_{k \sim 0} (1 - \eta(k)) \left( \frac{a\rho}{1-\rho} f(\sigma^k \eta) - b f(\eta) \right) + \eta(k) \left( \frac{b(1-\rho)}{\rho} f(\sigma^k \eta) - a f(\eta) \right),$$

where $\mathcal{L}_0$ is the same expression as $\mathcal{L}^{n,\rho}_{se}$ in (3.5) but the sum over $i \sim j$ is restricted to $i, j \in \Lambda_n \setminus \{0\}$.

Let $T^n_t(\nu_\rho)$ be the law at time $t$ of the process generated by either $\mathcal{L}^{n,\rho}_{se}$ or $\mathcal{L}^{n,\rho}_\beta$ conditioned on $\{\tau > t\}$ with initial measure $\nu_\rho$.

LEMMA 3.2. *Let $\{\underline{\nu}_n\}$ and $\{\overline{\nu}_n\}$ be two sequences of measures converging, respectively, to $\underline{\nu}$ and $\overline{\nu}$.*

(i) *Assume $\underline{\nu}_n \prec T^n_t(\nu_\rho) \prec \overline{\nu}_n$ for all $n$. Then $T^n_t(\nu_\rho)$ converges weakly to $T_t(\nu_\rho)$ and*

$$\underline{\nu} \prec T_t(\nu_\rho) \prec \overline{\nu}.$$

(ii) *Let $u_n$ be the unique positive principal eigenfunction of $(\mathcal{L}^{n,\rho}_{ab})^*$ with $\int u_n d\nu_\rho = 1$. Note that $(\mathcal{L}^{n,\rho}_{ab})^* u_n = 0$ and $d\mu_n := u_n d\nu_\rho$ is invariant for $\mathcal{L}^{n,\rho}_{ab}$. Assume that for all $n$, $\psi_n = d\underline{\nu}_n/d\nu_\rho$ is positive and decreasing (resp. $\psi'_n = d\overline{\nu}_n/d\nu_\rho$ is positive and increasing) such that $u_n/\psi_n$ is increasing (resp. $u_n/\psi'_n$ is decreasing). Assume also that $\underline{\nu}_n$ and $\overline{\nu}_n$ satisfy the FKG inequality. Then, there is a subsequence $\{n_k\}$ such that $d\mu_{n_k} := u_{n_k} d\nu_\rho$ converges weakly to $d\mu_\rho$, an invariant measure for $\mathcal{L}_{ab}$ with*

$$\underline{\nu} \prec \mu_\rho \prec \overline{\nu}.$$

PROOF. (i) We drop the subscripts $se$ or $\beta$ from the generators to unify their treatment. The stopped generator on $\mathcal{A}$, $\overline{\mathcal{L}}^{n,\rho} := \mathbb{1}_{\mathcal{A}^c} \mathcal{L}^{n,\rho}$ is bounded on $\Omega_n$ and it is obvious that

$$\forall \varphi \in \mathcal{D} \cap \{\varphi|_{\mathcal{A}} = 0\}, \ \forall \eta \in \Omega \quad \overline{\mathcal{L}}^{n,\rho} \varphi(\eta) \stackrel{n \to \infty}{\longrightarrow} \overline{\mathcal{L}} \varphi(\eta).$$

Thus, by a theorem of Trotter and Kurtz (see [9], Chapter I, Theorem 2.12), we have, for any $t \geq 0$,

$$(3.7) \quad P^{n,\rho}_\eta(\tau > t) = e^{t\overline{\mathcal{L}}^{n,\rho}} (\mathbb{1}_{\mathcal{A}^c})(\eta) \stackrel{n \to \infty}{\longrightarrow} e^{t\overline{\mathcal{L}}} (\mathbb{1}_{\mathcal{A}^c})(\eta) = P_\eta(\tau > t).$$

Note now that $T_t^n(\nu_\rho)$ is absolutely continuous with respect to $\nu_\rho$ and

$$(3.8) \qquad \frac{dT_t^n(\nu_\rho)}{d\nu_\rho}(\eta) = \frac{\exp(t\mathbb{1}_{\mathcal{A}^c}\mathcal{L}^{n,\rho})\mathbb{1}_{\mathcal{A}^c}(\eta)}{P_{\nu_\rho}^{n,\rho}(\tau > t)} = \frac{P_\eta^{n,\rho}(\tau > t)}{P_{\nu_\rho}^{n,\rho}(\tau > t)}.$$

Thus, by (3.7), (3.8) and dominated convergence, $T_t^n(\nu_\rho)$ converges weakly to $T_t(\nu_\rho)$ and point (i) follows easily.

(ii) Since $u_n/\psi_n$ is increasing, we have by the FKG inequality that for any increasing function $\varphi$,

$$(3.9) \qquad \int \varphi \, d\mu_n = \int \varphi \frac{u_n}{\psi_n} \, d\underline{\nu}_n \geq \int \varphi \, d\underline{\nu}_n \int \frac{u_n}{\psi_n} \, d\underline{\nu}_n = \int \varphi \, d\underline{\nu}_n,$$

so that $\mu_n \succ \underline{\nu}_n$. Similarly, we obtain that $\mu_n \prec \overline{\nu}_n$. Where as the space $\Omega^*$ is compact, there is a subsequence $\{n_k\}$ such that $\mu_{n_k}$ converges to a measure $\mu_\rho$. Now, for any function $\varphi \in \mathcal{D}$, $\mathcal{L}_{ab}^{n,\rho}\varphi$ converges to $\mathcal{L}_{ab}\varphi \in \mathcal{D}$. Thus, for $\varphi \in \mathcal{D}$,

$$0 = \int \mathcal{L}_{ab}^{n_k,\rho}(\varphi) \, d\mu_{n_k} \overset{k\to\infty}{\longmapsto} \int \mathcal{L}_{ab}(\varphi) \, d\mu_\rho.$$

Thus, $\int \mathcal{L}_{ab}(\varphi) \, d\mu_\rho = 0$ and $\mu_\rho$ is an invariant measure for $\mathcal{L}_{ab}$ with $\underline{\nu} \prec \mu_\rho \prec \overline{\nu}$.

$\square$

3.3. *Proof of Proposition* 1.3. The upper bound $T_t(\nu_\rho) \prec \nu_\rho$ is simple. Indeed, by observing that $\eta \mapsto P_\eta(\tau > t)$ is decreasing and by using the FKG inequality, we get, for any increasing $\varphi$,

$$\int \varphi \, dT_t(\nu_\rho) = \frac{1}{P_{\nu_\rho}(\tau > t)} \int \varphi(\eta) P_\eta(\tau > t) \, d\nu_\rho(\eta) \leq \int \varphi \, d\nu_\rho.$$

We now prove the lower bound $\nu_\alpha \prec T_t(\nu_\rho)$. First, notice that we are committing now an abuse of notation with $\prec$, since the monotonicity is only meant on $\mathcal{A}^c$. Henceforth, by $\mu \prec \nu$, for $\mu$ and $\nu$ with support in $\mathcal{A}^c$, we mean that for any $\varphi$ increasing on $\mathcal{A}^c$, $\int \varphi \, d\mu \prec \int \varphi \, d\nu$.

By Lemma 3.2, we need to establish two points: (a) for any integer $n$ and $t > 0$, $\underline{\nu}_n \prec T_t^n(\nu_\rho)$ and (b) that $\underline{\nu}_n$ tends to $\nu_\alpha$. Moreover, for (a), it is enough to show that

$$\eta \mapsto \frac{P_\eta^{n,\rho}(\tau > t)}{\psi_n(\eta)} \text{ is increasing on } \mathcal{A}^c \qquad \text{where } \psi_n = \frac{d\underline{\nu}_n}{d\nu_\rho}.$$

Indeed, note that on $\mathcal{A}^c$ the probability measure $\underline{\nu}_n$ satisfies Holley's condition (see [9], Theorem 2.9, page 75) which implies that $\underline{\nu}_n$ satisfies the FKG inequality. Thus, for any increasing function $\varphi$ on $\mathcal{A}^c$,

$$\int \varphi \, dT_t^n(\nu_\rho) = \int \varphi(\eta) \frac{P_\eta^{n,\rho}(\tau > t)}{\psi_n(\eta) P_{\nu_\rho}^{n,\rho}(\tau > t)} \, d\underline{\nu}_n(\eta)$$

$$\geq \int \varphi \, d\underline{\nu}_n \int \frac{P_\eta^{n,\rho}(\tau > t)}{\psi_n(\eta) P_{\nu_\rho}^{n,\rho}(\tau > t)} \, d\underline{\nu}_n(\eta) = \int \varphi \, d\underline{\nu}_n.$$

Now, we set

$$\psi_n(\eta) := \frac{1}{Z_n}(1 - \eta(0)) \prod_{i \in \Lambda_n \setminus \{0\}} \gamma_{i,n}^{\eta(i)}, \tag{3.10}$$

where $Z_n$ is a constant such that $\int \psi_n \, d\nu_\rho = 1$. Also, set

$$\alpha_i^{(n)} = \frac{\rho \gamma_{i,n}}{\rho \gamma_{i,n} + 1 - \rho}. \tag{3.11}$$

Note that (a) follows when the hypotheses of Lemma 2.1 are satisfied, whereas (b) follows as soon as for all sites $i$, $\alpha_i^{(n)} \to \alpha_i$, with $\sum (1 - \alpha_i/\rho)^2 < +\infty$.

We focus now on the four hypotheses of Lemma 2.1. Whereass the $\gamma_{i,n}$ are chosen smaller than 1, $\psi_n$ is decreasing. Moreover, a simple computation shows that $\mathcal{L}_{\psi_n}$, obtained by $\mathcal{L}_{se}^{n,\rho}$ as in Lemma 2.1, generates a monotone exclusion process since the intensity rate of any bond $(i,j)$ depends only on $\eta(i)$ and $\eta(j)$. Thus, it remains to show that $\mathcal{L}_{se}^{n,\rho}(\psi_n)/\psi_n$ is increasing.

Before specifying the $\{\gamma_i, i \in \mathbb{Z}^d\}$, we need some notation. Henceforth, we write $\mathcal{L}_n$ for $\mathcal{L}_{se}^{n,\rho}$ and $\gamma_i$ for $\gamma_{i,n}$. For each $i \in \mathbb{Z}^d$, let $\{X(i,t), t \geq 0\}$ be a symmetric simple random walk trajectory starting at $i$; we denote by $\mathbb{P}_i$ the average over such trajectory. Let

$$H_0 = \inf\{t : X(i,t) = 0\} \quad \text{and} \quad H_n = \inf\{t : X(i,t) \in \Lambda_n^c\}. \tag{3.12}$$

It is well known that for $i \sim 0$, $\mathbb{P}_i(H_0 < \infty) < 1/2$ for $d \geq 3$ (see, e.g., [5]) and that $\mathbb{P}_i(H_0 < H_n)$ increases to $\mathbb{P}_i(H_0 < \infty)$. Finally, note that $i \mapsto \mathbb{P}_i(H_0 < H_n)$ is harmonic outside 0 [see, e.g., (3.16)]. Let $0'$ be a neighbor of 0 and for $i \in \Lambda_n \setminus \{0\}$, set

$$\gamma_i = \frac{1}{1 + C_d \mathbb{P}_i(H_0 < H_n)} \quad \text{where} \quad C_d = \frac{1}{1 - 2\mathbb{P}_{0'}(H_0 < \infty)}. \tag{3.13}$$

Note that the corresponding $\alpha_i^{(n)}$—given through (3.11)—is

$$\alpha_i^{(n)} = \frac{\rho}{1 + (1-\rho)C_d \mathbb{P}_i(H_0 < H_n)}$$
$$\stackrel{n \to \infty}{\longrightarrow} \alpha_i := \frac{\rho}{1 + (1-\rho)C_d \mathbb{P}_i(H_0 < +\infty)}. \tag{3.14}$$

Thus, (b) follows as soon as $\sum_i \mathbb{P}_i^2(H_0 < +\infty) < +\infty$, that is, for $d \geq 5$ (see [2]).

PROOF THAT $V := \mathcal{L}_n \psi_n/\psi_n$ IS INCREASING. For $k \in \Lambda_n \setminus \{0\}$ and $\eta(k) = 0$ we show that $V(\sigma^k \eta) \geq V(\eta)$. We denote $\mathcal{N}_k := \{j \in \Lambda_n \setminus \{0\} : j \sim k\}$, $\mathcal{N}_k^0 := \{j \in \mathcal{N}_k : \eta(j) = 0\}$ and $\mathcal{N}_k^1 := \{j \in \Lambda_n \setminus \{0\} : \eta(j) = 1\}$. We treat the cases $k \in \Lambda_n \setminus \{\partial \Lambda_n, \mathcal{N}_0\}$, $k \in \partial \Lambda_n$ and $k \in \mathcal{N}_0$ separately.

*Case* 1. $k \in \Lambda_n \setminus \{\partial \Lambda_n, \mathcal{N}_0\}$. We assume $\eta(k) = 0$:

$$V(\sigma^k \eta) - V(\eta) = \sum_{j \sim k} \left( \left( \frac{\gamma_j}{\gamma_k} \right)^{1-\eta(j)} - 1 \right) - \sum_{j \sim k} \left( \left( \frac{\gamma_k}{\gamma_j} \right)^{\eta(j)} - 1 \right)$$

(3.15)

$$= \left( \sum_{j \in \mathcal{N}_k^0} \frac{\gamma_j}{\gamma_k} - |\mathcal{N}_k^0| \right) - \left( \sum_{j \in \mathcal{N}_k^1} \frac{\gamma_k}{\gamma_j} - |\mathcal{N}_k^1| \right).$$

Note that $i \mapsto 1/\gamma_i$ is harmonic at $k$, so that

(3.16) $$\sum_{j \in \mathcal{N}_k} \frac{1}{\gamma_j} = \frac{|\mathcal{N}_k|}{\gamma_k}.$$

Thus,

$$V(\sigma^k \eta) - V(\eta) = \sum_{j \in \mathcal{N}_k^0} \left( \frac{\gamma_j}{\gamma_k} + \frac{\gamma_k}{\gamma_j} \right) - 2|\mathcal{N}_k^0| \geq 0$$

(3.17)

$$\text{since for } x > 0, \ x + \frac{1}{x} \geq 2.$$

*Case* 2. $k \in \partial \Lambda_n$. Note that for any $\eta$,

$$\sqrt{\frac{d\sigma^k \nu_\rho}{d\nu_\rho}(\eta)} = \kappa^{2\eta(k)-1} \quad \text{with } \kappa := \sqrt{\left( \frac{1-\rho}{\rho} \right)} \quad \text{and}$$

(3.18)

$$\frac{\sigma^k \psi_n}{\psi_n} = \gamma_k^{1-2\eta(k)}.$$

Thus, for $\eta$ with $\eta(k) = 0$,

$$V(\sigma^k \eta) - V(\eta) = \sum_{j \in \mathcal{N}_k \cap \Lambda_n} \left( \left( \frac{\gamma_j}{\gamma_k} \right)^{1-\eta(j)} - \left( \frac{\gamma_k}{\gamma_j} \right)^{\eta(j)} \right)$$

(3.19)

$$+ n(k) \kappa \left( \frac{1}{\gamma_k} - 1 \right) - n(k) \frac{1}{\kappa} (\gamma_k - 1).$$

If we extend $\eta$ outside $\Lambda_n$ by 1 and recall that $\gamma_j = 1$ for $j \notin \Lambda_n$, we can replace the sum over $\mathcal{N}_k \cap \Lambda_n$ by a sum over $\mathcal{N}_k$ with an additional term $-n(k)(1 - \gamma_k)$. Thus,

$$V(\sigma^k \eta) - V(\eta) = \sum_{j \in \mathcal{N}_k} \left( \left( \frac{\gamma_j}{\gamma_k} \right)^{1-\eta(j)} - \left( \frac{\gamma_k}{\gamma_j} \right)^{\eta(j)} \right)$$

(3.20)

$$+ n(k)(1 - \gamma_k) \left( \frac{\kappa}{\gamma_k} + \frac{1}{\kappa} - 1 \right).$$

The same argument as in Case 1 implies that the sum over $\mathcal{N}_k$ is nonnegative and it is enough to have

$$(3.21) \quad \frac{\kappa}{\gamma_k} + \frac{1}{\kappa} - 1 \geq 0 \iff \frac{1}{\gamma_k} \geq \frac{1}{\kappa}\left(1 - \frac{1}{\kappa}\right),$$

which is always true for any $\rho \in ]0,1[$, since $\gamma_k \leq 1$.

*Case* 3. $k \in \mathcal{N}_0$. Note that for $\eta \notin \mathcal{A}$,

$$(3.22) \quad T^{k,0}\psi_n(\eta) = \begin{cases} 0, & \text{if } \eta(k) = 1, \\ \psi_\Lambda(\eta), & \text{if } \eta(k) = 0. \end{cases}$$

Thus, for $\eta(k) = 0$,

$$(3.23) \quad V(\sigma^k \eta) - V(\eta) = \left(\sum_{j \in \mathcal{N}_k^0} \frac{\gamma_j}{\gamma_k} - |\mathcal{N}_k^0|\right) - 1 - \left(\sum_{j \in \mathcal{N}_k^1} \frac{\gamma_k}{\gamma_j} - |\mathcal{N}_k^1|\right).$$

Now, whereas $i \mapsto \mathbb{P}_i(H_0 < H_n)$ is harmonic (and $0 \notin \mathcal{N}_k$ by definition),

$$(3.24) \quad \sum_{j \in \mathcal{N}_k} \frac{1}{\gamma_j} + (1 + C_d) = \frac{|\mathcal{N}_k| + 1}{\gamma_k}.$$

Thus, for our choice of $C_d$,

$$(3.25) \quad V(\sigma^k \eta) - V(\eta) \geq \frac{1 + C_d}{1 + C_d \mathbb{P}_k(H_0 < H_n)} - 2 \geq 0. \qquad \square$$

3.4. *Proof of Proposition* 1.4. Most of the arguments in the proof of Proposition 1.4 follow those in Section 3.3. A new difficulty arises from the fact that monotonicity of $\mathcal{L}_\psi$ is not trivial anymore.

As in Section 3.3, we first need some notations to specify the $\{\gamma_i, i \in \mathbb{Z}^d\}$. We denote by $\mathbb{P}_i$ the average over $\{X(i,t), t \geq 0\}$, a symmetric simple random walk trajectory starting at $i$. Let

$$(3.26) \, H_{\{0,0'\}} = \inf\{t: X(i,t) \in \{0,0'\}\} \quad \text{and} \quad H_n = \inf\{t: X(i,t) \in \Lambda_n^c\}.$$

We show in the Appendix that for $i \sim 0$, $i \neq 0'$, $\mathbb{P}_i(H_{\{0,0'\}} < \infty) < 1/2$ for $d \geq 4$. As $\mathbb{P}_i(H_{\{0,0'\}} < H_n)$ increases to $\mathbb{P}_i(H_{\{0,0'\}} < \infty)$, we choose $n$ large enough so that $\mathbb{P}_i(H_{\{0,0'\}} < H_n) < 1/2$. Thus,

$$C_2 := \sup_{k \in \mathcal{N}_0 \setminus \{0'\}} \frac{1}{1 - 2\mathbb{P}_k(H_{\{0,0'\}} < \infty)} > 0,$$

so that for all $k \in \mathcal{N}_0 \setminus \{0'\}$,

$$(3.27) \quad \frac{1 + C_2}{1 + C_2 \mathbb{P}_k(H_{\{0,0'\}} < H_n)} > 2.$$

We choose

$$\forall\, i \in \mathbb{Z}^d \qquad \gamma_i := \frac{1}{1 + C_2 \mathbb{P}_i(H_{\{0,0'\}} < H_n)} \quad \text{and}$$

(3.28)
$$\psi_n(\eta) = \mathbb{1}_{\mathcal{A}^c}(\eta) \prod_{i \in \Lambda_n} \gamma_i^{\eta(i)}.$$

Finally, note that $i \mapsto \mathbb{P}_i(H_{\{0,0'\}} < H_n)$ is harmonic outside $\{0, 0'\}$ and that with $\gamma_0 = \gamma_{0'}$, we have $T^b \psi_n = \psi_n$.

Define $\mathcal{L}_\beta^{n,\rho} = \mathcal{L}_{se}^{n,\rho} + (\beta - 1)(T^b - 1)$. We are now ready for the following proof.

PROOF THAT $V := \mathcal{L}_\beta^{n,\rho}(\psi_n)/\psi_n$ IS INCREASING. In the case where $k$ is not a neighbor of $0$ or of $0'$, then $V(\sigma^k \eta) - V(\eta)$ has the same expression as in Case 1 or 2 of Section 3.3. We do not repeat the computations.

*Case 1.* $k \in \mathcal{N}_0 \setminus \{0'\}$. Set $\mathcal{N}_k^* := \{j : j \sim k, j \notin \{0, 0'\}\}$, and for $\eta \notin \mathcal{A}$ and $\eta(k) = 0$,

(3.29)
$$V(\sigma^k \eta) - V(\eta) = S_k + \mathbb{1}_{\{\eta(0)=0, \eta(0')=0\}} \left(\frac{\gamma_0}{\gamma_k} - 1\right)$$
$$- \mathbb{1}_{\{\eta(0)=0, \eta(0')=1\}} - \mathbb{1}_{\{\eta(0)=1, \eta(0')=0\}} \left(\frac{\gamma_k}{\gamma_0} - 1\right)$$

with

(3.30)
$$S_k := \sum_{j \in \mathcal{N}_k^*} \left(\left(\frac{\gamma_j}{\gamma_k}\right)^{1-\eta(j)} - 1\right) - \sum_{j \in \mathcal{N}_k^*} \left(\left(\frac{\gamma_k}{\gamma_j}\right)^{1-\eta(j)} - 1\right).$$

Note that by harmonicity

(3.31)
$$\sum_{j \in \mathcal{N}_k^*} \frac{1}{\gamma_j} + \frac{1}{\gamma_{0'}} = \frac{|\mathcal{N}_k^*| + 1}{\gamma_k}.$$

Now, if we set $\mathcal{N}_k^0 := \{j \in \mathcal{N}_k^* : \eta(j) = 0\}$ and $\mathcal{N}_k^1 := \{j \in \mathcal{N}_k^* : \eta(j) = 1\}$, and use (3.31), the expression $S_k$ of (3.30) has the lower bound

(3.32)
$$S_k = \sum_{j \in \mathcal{N}_k^0} \frac{\gamma_j}{\gamma_k} - |\mathcal{N}_k^0| - \left(-\sum_{j \in \mathcal{N}_k^0} \frac{\gamma_k}{\gamma_j} + |\mathcal{N}_k^*| + 1 - \frac{\gamma_k}{\gamma_{0'}}\right)$$
$$= \sum_{j \in \mathcal{N}_k^0} \left(\frac{\gamma_j}{\gamma_k} + \frac{\gamma_k}{\gamma_j}\right) - 2|\mathcal{N}_k^0| + \frac{\gamma_k}{\gamma_{0'}} - 1 \geq \frac{\gamma_k}{\gamma_{0'}} - 1.$$

Now, in the event $\{\eta(0) = 1, \eta(0') = 0\}$, (3.29) and (3.32) yield

$$(3.33) \qquad V(\sigma^k \eta) - V(\eta) \geq \frac{\gamma_k}{\gamma_{0'}} - 1 - \left(\frac{\gamma_k}{\gamma_0} - 1\right) = 0.$$

In the event $\{\eta(0) = 0, \eta(0') = 1\}$, we have

$$V(\sigma^k \eta) - V(\eta) \geq \frac{\gamma_k}{\gamma_{0'}} - 2 \geq 0$$

(3.34)

$$\text{since } \frac{1 + C_2}{1 + C_2 \mathbb{P}_k(H_{\{0,0'\}} < H_n)} \geq 2 \qquad \text{[by (3.27)]}.$$

Finally, in $\{\eta(0) = 0, \eta(0') = 0\}$, we have

$$(3.35) \qquad V(\sigma^k \eta) - V(\eta) \geq \frac{\gamma_k}{\gamma_{0'}} - 1 + \frac{\gamma_0}{\gamma_k} - 1 \geq \frac{\gamma_k}{\gamma_{0'}} - 2 \geq 0.$$

*Case* 2. $k \in \{0, 0'\}$. We assume $k = 0$ and $\eta(0) = \eta(0') = 0$:

$$V(\sigma^k \eta) - V(\eta) = \sum_{j \in \mathcal{N}_0^*} \left(\left(\frac{\gamma_j}{\gamma_0}\right)^{1-\eta(j)} - 1\right) + (0 - |\mathcal{N}_{0'}^1|)$$

(3.36)
$$- \sum_{j \in \mathcal{N}_0^*} \left(\left(\frac{\gamma_0}{\gamma_j}\right)^{\eta(j)} - 1\right) - \sum_{j \in \mathcal{N}_{0'}^1} \left(\left(\frac{\gamma_0}{\gamma_j}\right)^{\eta(j)} - 1\right)$$

$$= \sum_{j \in \mathcal{N}_0^0} \frac{\gamma_j}{\gamma_0} - \sum_{j \in \mathcal{N}_0^1} \frac{\gamma_0}{\gamma_j} - \sum_{j \in \mathcal{N}_{0'}^1} \frac{\gamma_0}{\gamma_j} + |\mathcal{N}_0^1| - |\mathcal{N}_0^0|.$$

Condition (3.27) implies that $\gamma_j \geq 2\gamma_0$ for $j \in \mathcal{N}_0^* \cup \mathcal{N}_{0'}^*$. Thus,

$$(3.37) \qquad \begin{aligned} V(\sigma^k \eta) - V(\eta) &\geq 2|\mathcal{N}_0^0| - |\mathcal{N}_0^0| - (\tfrac{1}{2}|\mathcal{N}_0^1| - |\mathcal{N}_0^1|) - \tfrac{1}{2}|\mathcal{N}_{0'}^1| \\ &\geq |\mathcal{N}_0^0| + \tfrac{1}{2}|\mathcal{N}_0^1| - \tfrac{1}{2}(|\mathcal{N}_0^0| + |\mathcal{N}_0^1|) \geq 0. \end{aligned} \qquad \square$$

PROOF THAT $\mathcal{L}_{\psi_n}$ IS MONOTONE. We describe an order-preserving coupling between two trajectories $(\eta_t, \tilde{\eta}_t)$ for $t \geq 0$, when $\eta_0 \succ \tilde{\eta}_0$. We run the two dynamics with the same family of Poisson processes up to the first time there is a mismatch at $0$ or $0'$. Assume that this happens at the stopping time $T$ and that $\eta_T(0) = 1 = 1 - \tilde{\eta}_T(0)$. Under $\mathcal{L}_{\psi_n}$, the rate for bringing $\eta$ particles from any site of $\mathcal{N}_{0'}$ to $0'$ is null. Let $\{\tilde{\tau}_i \circ \theta_T, i \in \mathcal{N}_{0'} \setminus \{0\}\}$ be the exponential times associated with the bonds of $0'$ in $\tilde{\eta}$ after time $T$. Note that if $\tilde{\eta}_T(i) = 1$ for the $i$ neighbor of $0'$, then the intensity rate of $(i, 0')$ is $\gamma_0/\gamma_i \leq 1$. Thus,

$$\alpha := \sum_{\substack{i \sim 0' \\ i \neq 0}} \tilde{\eta}_T(i) \frac{\gamma_0}{\gamma_i} \leq 2d - 1.$$

Let $\tau_a$ be an exponential time of parameter $\beta - \alpha$ independent of the other times. We associate to the bond $b$ of $\eta$ at time $T$ the exponential time of parameter $\beta$:

$$(3.38) \qquad \tau_b := \min(\tau_a, \{\tilde{\tau}_i \circ \theta_T, i \in \mathcal{N}_{0'} \setminus \{0\}, \tilde{\eta}_T(i) = 1\}).$$

We associate to the bond $b$ of $\tilde{\eta}$ an independent copy of $\tau_b$, but since $\tilde{\eta}_T(0) = \tilde{\eta}_T(0') = 0$, this has no effect. All remaining bonds in the two trajectories share the same Poisson processes. Now, if $\tau_b = \tau_a$, then there is a mismatch at 0 and $\eta_{T+\tau_b+} \succ \tilde{\eta}_{T+\tau_b+}$, and we restart the same construction, with 0 and $0'$ exchanging roles. On the other hand, if $\tau_b < \tau_a$, then the mismatch at 0 and $0'$ vanishes, $\eta_{T+\tau_b+} \succ \tilde{\eta}_{T+\tau_b+}$, and we proceed with the same Poisson processes on all bonds. □

3.5. *Proof of Proposition* 1.6. We rely here on Lemmas 2.2 and 3.2, with $\mathcal{L} = \mathcal{L}_{ab}^{n,\rho}$. We define

$$(3.39) \qquad \psi_n(\eta) = \frac{1}{Z_n} \prod_{i \in \Lambda_n \setminus \{0\}} \gamma_i^{\eta(i)},$$

where $Z_n$ is a constant such that $\int \psi_n \, d\nu_\rho = 1$ and

$$(3.40) \qquad \gamma_i = \frac{1}{1 + C_{a,b} \mathbb{P}_i(H_0 < H_n)},$$

where $H_0$ and $H_n$ are defined in (3.12), and $C_{a,b}$ is a constant that will be fixed later. We remark that

$$(\mathcal{L}_{ab}^{n,\rho})^* f = \tilde{\mathcal{L}} f + \left(\frac{a}{1-\rho} - \frac{b}{\rho}\right) \sum_{k \sim 0} (\rho - \eta(k)) f,$$

where $\tilde{\mathcal{L}}$ is the Markov generator

$$(3.41) \quad \tilde{\mathcal{L}} f = \mathcal{L}_0 f + \sum_{k \sim 0} (1 - \eta(k)) \frac{a\rho}{1-\rho}(\sigma^k f - f) + \eta(k) \frac{b(1-\rho)}{\rho}(\sigma^k f - f).$$

Thus, as observed in [6], if $a\rho = b(1-\rho)$, then $(\mathcal{L}_{ab}^{n,\rho})^*$ is a Markov generator and $\nu_\rho$ is an invariant measure (reversible if $a = b$).

Since $\psi$ is a product function, $\mathcal{L}_\psi$ is a monotone generator. Indeed, the intensity rate of $(i,j)$ depends only on $\eta(i)$ and $\eta(j)$, whereas the rate of spin flip at site $k$ depends only on $\eta(k)$. Thus, to prove the lower bound in (1.7), we are left to show the following proof.

PROOF THAT $V := (\mathcal{L}_{ab}^{n,\rho})^* \psi_n / \psi_n$ IS INCREASING. We take $k \in \Lambda_n \setminus \{0\}$ with $\eta(k) = 0$ and we show that $V(\sigma^k \eta) - V(\eta) \geq 0$. The case where $k$ is not a neighbor of 0 is similar to Cases 1 or 2 in the proof of Proposition 1.3. Assume $k \sim 0$. Rewriting $V$, we need

$$
\begin{aligned}
\sum_{j \in \mathcal{N}_k^0} \frac{\gamma_j}{\gamma_k} + |\mathcal{N}_k^1| + \left(\frac{b(1-\rho)}{\rho}\frac{1}{\gamma_k} - a\right) \\
\geq \sum_{j \in \mathcal{N}_k^1} \frac{\gamma_k}{\gamma_j} + |\mathcal{N}_k^0| + \left(\frac{a\rho}{1-\rho}\gamma_k - b\right).
\end{aligned}
\tag{3.42}
$$

By defining $\gamma_0 = 1/(1 + C_{a,b})$, we obtain that $k \mapsto 1/\gamma_k$ is harmonic outside 0 and we obtain the sufficient condition

$$
b\left(\frac{1-\rho}{\rho}\frac{1}{\gamma_k} + 1\right) - a\left(1 + \frac{\rho}{1-\rho}\gamma_k\right) \geq 1 - \frac{\gamma_k}{\gamma_0}.
\tag{3.43}
$$

For the upper bound in (1.7), we replace $\psi_n$ with

$$
\psi_n'(\eta) = \frac{1}{Z_n'} \prod_{i \in \Lambda_n \setminus \{0\}} \gamma_i^{-\eta(i)},
\tag{3.44}
$$

where $Z_n'$ is a constant such that $\int \psi_n' \, d\nu_\rho = 1$. It is easy to check that the corresponding $\tilde{\alpha}_i = (\rho \gamma_i^{-1})/(\rho \gamma_i^{-1} + 1 - \rho)$ produces the relationship in (1.7). In this case, the corresponding potential $V'$ should be decreasing. By the same argument used above, we obtain the sufficient condition

$$
a\left(1 + \frac{\rho}{1-\rho}\frac{1}{\gamma_k}\right) - b\left(\frac{1-\rho}{\rho}\gamma_k + 1\right) \geq 1 - \frac{\gamma_k}{\gamma_0}.
\tag{3.45}
$$

Now, if we set $\delta = b(1-\rho)/(a\rho)$, (3.43) and (3.45) read

$$
\left(\frac{\delta}{\gamma_k} - 1\right)\left(a + \frac{\gamma_k}{\delta}b\right) \geq -\frac{C_{a,b}\mathbb{P}_k(H_0 > H_n)}{1 + C_{a,b}\mathbb{P}_k(H_0 < H_n)}
\tag{3.46}
$$

and

$$
\left(\frac{1}{\gamma_k \delta} - 1\right)(b + \gamma_k \delta a) \geq -\frac{C_{a,b}\mathbb{P}_k(H_0 > H_n)}{1 + C_{a,b}\mathbb{P}_k(H_0 < H_n)}.
\tag{3.47}
$$

Thus, for any $a$ and $b$ positive, we can take $C_{a,b}$ large enough so that (3.46) and (3.47) hold. $\square$

**4. Proof of Proposition 1.8.** Define $\overline{S}_t = \mathbb{1}_{\mathcal{A}^c} \exp[t \mathbb{1}_{\mathcal{A}^c} \mathcal{L}]$. Let us first note that (ii) is a simple consequence of (i). Indeed, let $g = d\pi/d\nu$ and let $f$ be in $L^2(\nu)$:

$$
\begin{aligned}
\frac{\int \overline{S}_t f \, d\pi}{P_\pi(\tau > t)} &= \frac{\int \overline{S}_t f \, d\pi}{P_\nu(\tau > t)} \frac{P_\nu(\tau > t)}{\int \overline{S}_t g \, d\nu} \\
&= \frac{E_\nu[g(\eta_0) f(\eta_t) \mathbb{1}_{\tau > t}]}{P_\nu(\tau > t)} \left(\frac{E_\nu[g(\eta_t) \mathbb{1}_{\tau > t}]}{P_\nu(\tau > t)}\right)^{-1} \xrightarrow{t \to \infty} \frac{\int f \, d\mu \int g \, d\mu}{\int g \, d\mu}.
\end{aligned}
$$

Now, to prove (i), we first set

$$H_t = \frac{\overline{S}_t g}{P_\nu(\tau > t)} \quad \text{and} \quad H = u \int g \, d\mu,$$

and we need to show that $H_t$ converges to $H$ in the weak-$L^2(\nu)$ topology. We actually show that this convergence holds in $L^2(\nu)$, which is equivalent to the two facts

$$\lim_{t \to \infty} \int H_t H \, d\nu = \int H^2 \, d\nu \tag{4.1}$$

and

$$\lim_{t \to \infty} \int H_t^2 \, d\nu = \int H^2 \, d\nu. \tag{4.2}$$

We begin by proving (4.1). Since $u$ is a simple eigenfunction in $L^2(\nu)$, $\overline{S}_t(u) = e^{-\lambda t} u$ $\nu$-a.s. and, by symmetry,

$$\int H_t H \, d\nu = \int u \overline{S}_t(g) \, d\nu \frac{\int g \, d\mu}{P_\nu(\tau > t)} = \int g \overline{S}_t(u) \, d\nu \frac{\int g \, d\mu}{P_\nu(\tau > t)}$$
$$= \frac{e^{-\lambda t}}{P_\nu(\tau > t)} \left( \int g \, d\mu \right)^2 \overset{t \to \infty}{\longrightarrow} \int u^2 \, d\nu \left( \int g \, d\mu \right)^2 = \int H^2 \, d\nu.$$

In the last step, we used (1.4). Thus, (4.1) is established.

In order now to prove (4.2), we rewrite

$$\int H_t^2 \, d\nu = \frac{\int g \overline{S}_{2t} g \, d\nu}{P_\nu(\tau > t)^2} = \frac{\int g \overline{S}_{2t} g \, d\nu}{(\int \overline{S}_t g \, d\nu)^2} \frac{(\int \overline{S}_t g \, d\nu)^2}{P_\nu(\tau > t)^2}.$$

Since

$$\lim_{t \to \infty} \frac{\int \overline{S}_t g \, d\nu}{P_\nu(\tau > t)} = \int g \, d\mu,$$

we are left to show that

$$\lim_{t \to \infty} \frac{\int g \overline{S}_{2t} g \, d\nu}{(\int \overline{S}_t g \, d\nu)^2} = \int u^2 \, d\nu. \tag{4.3}$$

Denote by $(\Pi_x)_{x \in \mathbb{R}}$ the spectral projections of $\overline{\mathcal{L}}$ in $L^2_{\mathcal{A}}$. We know that $\Pi_x = I$ for $x \geq -\lambda$. Thus, by the spectral theorem,

$$\int g \overline{S}_{2t} g \, d\nu = \int_{(-\infty, -\lambda]} e^{2tx} \, d\langle g, \Pi_x g \rangle, \tag{4.4}$$

where $\langle \cdot, \cdot \rangle$ in the scalar product in $L^2(\nu)$. Now, we have the orthogonal decomposition

$$g = \langle g, \bar{u} \rangle \bar{u} + \varphi \quad \text{with } \bar{u} = \frac{u}{\|u\|_2} \ (\|\cdot\|_2 = \|\cdot\|_{L^2(\nu)}).$$

By assumption $(\mathcal{H})$, $\lambda$ is a simple eigenvalue for $\overline{\mathcal{L}}$. This implies that $\text{range}(\Pi_{-\lambda} - \Pi_{-\lambda^-}) = \text{span}(u)$. Indeed, since the spectrum of $\overline{\mathcal{L}}$ is bounded from above, we have that $\text{range}(\Pi_{-\lambda} - \Pi_{-\lambda^-}) \subset D(\overline{\mathcal{L}})$, so that Theorem 5 on page 265 of [7] applies and $\varphi = \Pi_{-\lambda^-}(\varphi)$. In particular, $\varphi_n := \Pi_{-\lambda-1/n}\varphi$, converges to $\varphi$ in $L^2(\nu)$. Define
$$g_n = \langle g, \bar{u}\rangle \bar{u} + \varphi_n.$$
Since $\langle g_n, \Pi_x g_n\rangle = \langle \varphi_n, \Pi_x \varphi_n\rangle$ for $x < -\lambda$ and
$$\langle g_n, \Pi_{-\lambda} g_n\rangle - \langle g_n, \Pi_{-\lambda^-} g_n\rangle = \langle g_n, g_n\rangle - \langle \varphi_n, \varphi_n\rangle = \langle g, \bar{u}\rangle^2,$$
we have

(4.5)
$$\int_{(-\infty,-\lambda]} e^{2tx} \, d\langle g_n, \Pi_x g_n\rangle - e^{-2t\lambda}\langle g, \bar{u}\rangle^2$$
$$= \int_{(-\infty,-\lambda-1/n]} e^{2tx} \, d\langle \varphi_n, \Pi_x \varphi_n\rangle = o(e^{-2t\lambda}).$$

Similarly,

(4.6)
$$\int \overline{S}_t g_n \, d\nu = \int_{(\infty,-\lambda]} e^{tx} \, d\langle \mathbb{1}_{A^c}, \Pi_x g_n\rangle$$
$$= \int_{(\infty,-\lambda-1/n]} e^{tx} \, d\langle \mathbb{1}_{A^c}, \Pi_x \varphi_n\rangle + e^{-\lambda t}\langle g, \bar{u}\rangle\langle \bar{u}, \mathbb{1}_{A^c}\rangle$$
$$= e^{-\lambda t}\frac{\langle g, \bar{u}\rangle}{\|u\|_2} + o(e^{-\lambda t}).$$

By (4.4), (4.5) and (4.6), we have that (4.3) holds if we replace $g$ with $g_n$ and therefore (1.9) holds for $g_n$. To complete the proof, we are left to show that

(4.7)
$$\lim_{n\to\infty} \sup_t \left|\frac{\int f\overline{S}_t g \, d\nu}{P_\nu(\tau > t)} - \frac{\int f\overline{S}_t g_n \, d\nu}{P_\nu(\tau > t)}\right| = 0.$$

However,
$$\left|\frac{\int f\overline{S}_t g \, d\nu}{P_\nu(\tau > t)} - \frac{\int f\overline{S}_t g_n \, d\nu}{P_\nu(\tau > t)}\right| \leq \frac{\int |(g_n - g)\overline{S}_t f| \, d\nu}{P_\nu(\tau > t)}$$
$$\leq \frac{1}{P_\nu(\tau > t)}\|\overline{S}_t f\|_2 \|g_n - g\|_2$$
$$\leq \frac{e^{-\lambda t}}{P_\nu(\tau > t)}\|f\|_2 \|g_n - g\|_2.$$

The proof is concluded after recalling that
$$\|g_n - g\|_2 \to 0 \quad \text{and} \quad \sup_t \frac{e^{-\lambda t}}{P_\nu(\tau > t)} < \infty \qquad \text{(by Fact 1.2).}$$

**5. The process $P^u$.** In this section, we study the law of the whole path $\eta_{[0,t]} \equiv (\eta_s)_{s \in [0,t]}$ under the conditional distribution $P_\nu(\cdot|\tau > t)$, in the limit as $t$ tends to infinity. Consider the stochastic process

$$Z_t = \frac{u(\eta_0)u(\eta_t)e^{\lambda t}}{\int u^2\,d\nu}\mathbb{1}_{\tau>t}.$$

Let $\mathcal{F}_t$ be the $\sigma$-field $\sigma\{\eta_s : s \in [0,t]\}$. Note that, for $0 \le s < t$,

$$E_\nu(Z_t|\mathcal{F}_s) = \frac{u(\eta_0)e^{\lambda t}\mathbb{1}_{\tau>s}}{\int u^2\,d\nu}e^{(t-s)\overline{\mathcal{L}}}u(\eta_s) = Z_s, \qquad \nu\text{-a.s.},$$

so that $(Z_t)_{t\ge 0}$ is a positive martingale under $P_\nu$ with $E_\nu[Z_t] = 1$ for any $t \ge 0$. Thus, for any $t \ge 0$, a probability measure $P^u$ can be defined on $\mathcal{F}_t$ by

$$\left.\frac{dP^u}{dP_\nu}\right|_{\mathcal{F}_t} = Z_t.$$

Let $d\hat{\mu} = u^2\,d\nu/\int u^2\,d\nu$. For $g \in L^2(\hat{\mu})$ and $t \ge s \ge 0$, we have, using reversibility,

$$E^u[g(\eta_t)] = \int g(\eta_t)Z_t\,dP_\nu = \int \frac{e^{\lambda t}u\overline{S}_t(ug)}{\int u^2\,d\nu}\,d\nu = \int g\,d\hat{\mu}$$

and

(5.1)
$$\begin{aligned}E^u[g(\eta_t)|\mathcal{F}_s] &= \frac{1}{Z_s}E_\nu[Z_tg(\eta_t)|\mathcal{F}_s]\\ &= \frac{1}{Z_s}\frac{e^{\lambda t}u(\eta_0)\mathbb{1}_{\tau>s}}{\int u^2\,d\nu}\overline{S}_{t-s}(ug)(\eta_s) = \frac{e^{\lambda(t-s)}\overline{S}_{t-s}(ug)(\eta_s)}{u(\eta_s)},\end{aligned}$$

where equalities are intended $P^u$-a.s. Therefore, under $P^u$, the canonical process $\eta_t$ is stationary with marginal law $\hat{\mu}$ and the transition probabilities are given by

(5.2) $$\forall \xi \in \mathcal{A}^c \qquad E_\xi^u[g(\eta_t)] = \frac{e^{\lambda t}}{u(\xi)}\overline{S}_t(gu)(\xi).$$

By the same argument in (5.1), the associated Markov family $\{P_\xi^u, \xi \in \Omega\}$ is given on $\mathcal{A}^c$ by

(5.3) $$\forall \xi \in \mathcal{A}^c \qquad P_\xi^u(\{\eta : \eta_{[0,t]} \in \Gamma\}) = \frac{1}{u(\xi)}E_\xi[\mathbb{1}_\Gamma(\eta_{[0,t]})u(\eta_t)e^{\lambda t}\mathbb{1}_{\{\tau>t\}}],$$

where $\Gamma$ is a measurable set of paths depending only on times in $[0,t]$. Observe, finally, that $P^u$ is reversible, that is, it is invariant by time reversal.

5.1. *Proof of Proposition* 1.9. Let $\varphi = \varphi(\eta_{[0,r]})$ be a bounded measurable function. By reversibility and the Markov property,

$$E_\nu(\varphi(\eta_{[a_t, a_t+r]})|\tau > t)$$
$$= \frac{E_\nu(\varphi(\eta_{[a_t, a_t+r]})\mathbb{1}_{\{\tau>t\}}))}{P_\nu(\tau > t)}$$
$$= \frac{E_\nu(P_{\eta_0}(\tau > a_t)\varphi(\eta_{[0,r]})\mathbb{1}_{\{\tau>r\}}P_{\eta_r}(\tau > t - a_t - r))}{P_\nu(\tau > t)}$$
$$= E_\nu\left[\frac{dT_{a_t}(\nu)}{d\nu}(\eta_0)\frac{dT_{t-a_t-r}(\nu)}{d\nu}(\eta_r)\varphi(\eta_{[0,r]})\mathbb{1}_{\{\tau>r\}}\right]\beta(t)$$

with

$$\beta(t) = \frac{P_\nu(\tau > a_t)P_\nu(\tau > t - a_t - r)}{P_\nu(\tau > t)}$$
$$= e^{\lambda r}\frac{P_\nu(\tau > a_t)}{e^{-\lambda a_t}}\frac{e^{-\lambda t}}{P_\nu(\tau > t)}\frac{P_\nu(\tau > t - a_t - r)}{e^{-\lambda t - a_t - r}}.$$

Now, recalling (1.4),

(5.4) $$\beta(t) \stackrel{t\to\infty}{\longrightarrow} \frac{e^{-\lambda r}}{\int u^2\, d\nu}.$$

Also, by the Cauchy–Schwarz inequality, if we set $f(t, \eta) = (dT_t(\nu)/d\nu)(\eta)$, then

(5.5)
$$\left|\int\left(\frac{dT_{a_t}(\nu)}{d\nu}(\eta_0)\frac{dT_{t-a_t-r}(\nu)}{d\nu}(\eta_r) - u(\eta_0)u(\eta_r)\right)\varphi(\eta_{[0,r]})\mathbb{1}_{\{\tau>r\}}\, dP_\nu\right|$$
$$\leq |\varphi|_\infty\left(\int |f(a_t, \eta_0) - u(\eta_0)|u(\eta_r)\, dP_\nu\right.$$
$$\left.+ \int |f(t - a_t - r, \eta_r) - u(\eta_r)|u(\eta_0)\, dP_\nu\right)$$
$$\leq |\varphi|_\infty\left(\left(\int u^2(\eta_r)\, dP_\nu \int |f(a_t, \eta_0) - u(\eta_0)|^2\, dP_\nu\right)^{1/2}\right.$$
$$\left.+ \left(\int u^2(\eta_0)\, dP_\nu \int |f(t - a_t - r, \eta_r) - u(\eta_r)|^2\, dP_\nu\right)^{1/2}\right)$$
$$\leq |\varphi|_\infty\|u\|_2(\|f(a_t, \cdot) - u\|_2 + \|f(t - a_t - r, \cdot) - u\|_2).$$

This last expression goes to 0 as $t$ tends to infinity. Thus, gathering (5.4) and (5.5), we obtain

$$E_\nu(\varphi(\eta_{[a_t, a_t+r]})|\tau > t)$$

$$\stackrel{t\to\infty}{\longrightarrow} E_\nu[u(\eta_0)u(\eta_r)\mathbb{1}_{\{\tau>r\}}\varphi(\eta_{[0,r]})]\frac{e^{\lambda r}}{\int u^2\, d\nu} = E^u(\varphi(\eta_{[0,r]})).$$

REMARK 5.1. By using arguments as those in Section 4, we can show that, for

$$0 < a_t < b_t < t \quad \text{with } \lim_{t\to\infty} a_t = \lim_{t\to\infty}(b_t - a_t) = \lim_{t\to\infty}(t - b_t) = \infty,$$

the paths $\eta_{[a_t,a_t+r]}$ and $\eta_{[b_t,b_t+r]}$ decouple with respect to $P_\nu(\cdot|\tau>t)$ as $t\to\infty$, that is,

$$\lim_{t\to\infty} E_\nu[\varphi(\eta_{[a_t,a_t+r]})\psi(\eta_{[b_t,b_t+r]})|\tau>t] = E^u(\varphi(\eta_{[0,r]}))E^u(\psi(\eta_{[0,r]}))$$

for $\varphi,\psi$ bounded and measurable. In particular, the following generalization of (1.9) holds:

$$(5.6) \qquad \lim_{t\to\infty} \frac{E_\nu[f(\eta_{a_t})g(\eta_{b_t})\mathbb{1}_{\{\tau>t\}}]}{P_\nu(\tau>t)} = \int f\, d\hat\mu \int g\, d\hat\mu.$$

REMARK 5.2. Concerning the asymptotics at the boundary of $[0,t]$, we have the following result. For $r > 0$, the distribution of $\{\eta_s, s \in [0,r]\}$ with respect to $P_\nu(\cdot|\tau>t)$ converges to the restriction to the time interval $[0,r]$ of $P^u_\mu \equiv \int P^u_\xi\, d\mu$, while the distribution of $\{\eta_s, s \in [t-r,t]\}$ with respect to $P_\nu(\cdot|\tau>t)$ converges to the time reversal of the restriction to the time interval $[0,r]$ of $P^u_\mu \equiv \int P^u_\xi\, d\mu$. Indeed, by reversibility, the two statements above are equivalent, so we prove only the first one. The argument is identical to that in Proposition 1.9. For $\varphi = \varphi(\eta_{[0,r]})$ bounded and measurable, we have

$$E_\nu[\varphi(\eta_{[0,r]})|\tau>t] = \frac{E_\nu[\varphi(\eta_{[0,r]})\mathbb{1}_{\{\tau>r\}}P_{\eta_r}(\tau>t-r)]}{P_\nu(\tau>t)}$$

$$= E_\nu\left[\varphi(\eta_{[0,r]})\mathbb{1}_{\{\tau>r\}}\frac{dT_{t-r}(\nu)}{d\nu}(\eta_r)\right]\frac{P_\nu(\tau>t-r)}{P_\nu(\tau>t)}$$

$$\stackrel{t\to\infty}{\longrightarrow} E_\nu[\varphi(\eta_{[0,r]})\mathbb{1}_{\{\tau>r\}}u(\eta_r)e^{\lambda r}] = E^u_\mu[\varphi(\eta_{[0,r]})].$$

## APPENDIX

We show in this appendix that, with the notation of Section 3.4, if $k$ is a neighbor of $0$, $k \neq 0'$, then in dimensions $d \geq 4$,

$$\mathbb{P}_k(H_{\{0,0'\}} < \infty) < \tfrac{1}{2}, \qquad \text{where } H_\Lambda = \inf\{n>0 : S_n \in \Lambda\},$$

where $\Lambda \subset \mathbb{Z}^d$ and $\{S_n, n \in \mathbb{N}\}$ is a random walk. First, note that

$$\mathbb{P}_k(H_{\{0,0'\}} < \infty) \leq \mathbb{P}_k(H_0 < \infty) + \mathbb{P}_k(H_{0'} < \infty).$$

We will show that (i) $\mathbb{P}_k(H_{0'} < \infty) \leq \mathbb{P}_k(H_0 < \infty)$ and that (ii) $\mathbb{P}_k(H_0 < \infty) \leq \mathbb{P}_0(H_0 < \infty)$. Assume (i) and (ii) hold. If $R$ is the number of returns to the origin, we have the classical equality

$$\mathbb{P}_0(H_0 < \infty) = \frac{\mathbb{E}_0[R]}{1 + \mathbb{E}_0[R]} \qquad \left(\text{where we recall that } \mathbb{E}_0[R] = \sum_{n=2}^{\infty} \mathbb{P}_0(S_n = 0)\right).$$

Finally, we conclude, using the computation in [8], that $\mathbb{E}_0[R] < 0.25$ for $d \geq 4$.

Now, we show (i). To each path starting from $k$ and touching $0'$, we associate a path starting from $k$ and touching $0$. Let $\{S_n, n \in \mathbb{N}\}$ be a path with $S_0 = k$, let

$$\nu = \inf\{n > 0 : S_n - S_{n-1} = \overrightarrow{00'}\}$$

and note that $H_{0'} > \nu$. Define $\{S'_n, n \in \mathbb{N}\}$ as follows: if $\nu = \infty$, then $S'_n = S_n$ for all $n$; otherwise, let $S'_n = S_n$ for $n < \nu$ and $S'_n = S_{n+1} - \overrightarrow{00'}$ for $n \geq \nu$. Let $H'_0 = \inf\{n : S'_n = 0\}$. Note that if $H_{0'} < \infty$, then $S_{\nu-1} = 0$. Thus, $(S_n, S'_n)$ is a coupling where $H'_0 \leq H_{0'}$, and where each marginal is a random walk. Thus, (i) holds.

Now, point (ii). We couple $S_n$ with a path $\tilde{S}_n$ starting at 0 and such that if $S_n = 0$, then $\tilde{S}_{n+1} = 0$. For $i, j$ two sites that are neighbors of 0, let $R_{i,j}$ be the rotation with center 0 which sends $\overrightarrow{0i}$ onto $\overrightarrow{0j}$. Let $X_0$ be a uniform choice of a site in $\mathcal{N}_0$, and define

$$\tilde{S}_1 = X_0 \quad \text{and} \quad \text{for } n \geq 1, \qquad \tilde{S}_{n+1} = X_0 + R_{k,X_0}(S_n).$$

This definition ensures that $\{\tilde{S}_n, n \in \mathbb{N}\}$ has independent increments uniformly in $\mathcal{N}_0$ and such that if $S_n = 0$, then $\tilde{S}_{n+1} = 0$. Thus, (ii) follows easily.

CENTRE DE MATHEMATIQUES
ET INFORMATIQUE
UNIVERSITÉ DE PROVENCE
39 RUE JOLIOT CURIE
F-13453 MARSEILLE CEDEX 13
FRANCE
E-MAIL: asselah@cmi.univ-mrs.fr

DIPARTIMENTO DI MATEMATICA
PURA E APPLICATA
UNIVERSITÀ DI PADOVA
VIA BELZONI 7
35131 PADOVA
ITALY
E-MAIL: daipra@math.unipd.it